\documentclass{amsart}

\usepackage{mathdots,yhmath}

\newtheorem{thm}{Theorem}[section]

\newtheorem{cor}[thm]{Corollary}
\newtheorem{lem}[thm]{Lemma}

\theoremstyle{definition}

\theoremstyle{remark}

\numberwithin{equation}{section}

\def\qq{\qquad}

\def\0{\mathbf{0}}

\begin{document}

\title[The Generalized Remainder and Quotient Theorems]{The Generalized Remainder and Quotient Theorems of Univariate Polynomials}

\author{Wiwat Wanicharpichat}

\thanks{Received September 9, 2011.}
\thanks{ {\it Key  words and phrases.} Generalized remainder theorem, generalized quotient theorem, lower Hessenberg-Toeplitz matrix,
Hankel matrix,  recurrent sequence.
}
\thanks{Research supported by Naresuan University.}

%%%%%%%%%%%%%%%%%

\begin{abstract}
The author in \cite{WW10} was proved the generalized remainder and quotient theorems of polynomial in one indeterminate where the divisor is complete factorization to linear factors. In this paper we give the formula for the generalized remainder theorem and the generalized quotient theorem of polynomials  when the divisor is not factorization to linear factors.
\end{abstract}

\maketitle
\vspace{-9mm}
%%%%%%%%%%%%%%%%%%%%%%%%%%%%%%%%%%%%%%%%%%%%%%%%%%%%%%%%%%%%%%%%%%%

\section{Introduction and Preliminary}

Let $f(x)$ and $g(x)$ be polynomials over a field $F$. Since $F[x]$ is a Euclidean domain, the Division Algorithm holds and hence there are (unique) polynomials $q(x)$ and $r(x)$ in $F[x]$\ such that
\[
f(x)=g(x)q(x)+r(x)
\]
where $r(x)=0$ or its degree is less than that of $g(x).$ If $g(x) = x - b,$ the well known Remainder Theorem gives $r(x) = f(b).$  The author in \cite{WW85} was shown that, if $g(x) = ax - b,~a \ne 0,$ then a simple extension to this result gives the remainder as $f(a^{-1}b).$ This is true since, $f(x) = (x - a^{-1}b)q(x) + f(a^{-1}b),$ by the remainder theorem, for some polynomial $q(x)$ and by taking out a constant $a$ and writing $q(x)= aq_{1}(x),$ we have $f(x) = (ax - b)q_{1}(x) + f(a^{-1}b).$ Thus, the remainder term does not change and it is sufficient to consider divisors which are monic polynomials. This result may be extended to divisors which are polynomials of degree higher than $1.$ If $f(x)$ is a polynomial of degree $n$ and $g(x)$ is a polynomial of degree $m,$ the remainder polynomial $r(x)$ is unchanged if $g(x)$ is replaced by the corresponding monic polynomial obtained from it by taking out the leading coefficient.

\begin{lem}[{\rm\cite[Lemma 2.1]{WW10}}]\label{lem-S1.1}
Let $f(x)$\ and $g(x)= b_{m}x^{m}+b_{m-1}x^{m-1}+\dots+b_{1}x+ b_{0}$ be polynomials in $F[x]$, $\deg g(x)=m$ and $b_{m}\neq 1$, if
$$g_{1}(x)=\frac{1}{b_{m}}g(x)=  x^{m}+\beta_{m-1}x^{m-1}+\dots+\beta_{1}x+ \beta_{0}$$
where $\beta_{i}=\frac{b_{i}}{b_{m}}$, $i=0,1,\ldots, m-1$, be the corresponding monic polynomial of $g(x)$, and $q_{1}(x)$, $ q(x)$ be the quotients and $r_{1}(x),$\ $r(x)$ be the remainders on
dividing $f(x)$\ by $g_{1}(x)$,  and by $g(x)$ respectively, then
\begin{equation}
q(x)=\frac{1}{b_{m}}q_{1}(x)\quad\text{and}\quad r(x)=r_{1}(x).
\end{equation}
\end{lem}

\begin{proof} Since
\[
g(x)=b_{0}+b_{1}x+\ldots +b_{m-1}x^{m-1}+b_{m}x^{m},\quad
b_{m}\neq  1.
\]
Let
\[
\begin{array}{lcl}
 g_{1}(x)&=&\frac{1}{b_{m}}g(x)\\
        &=&\frac{b_{0}}{b_{m}}+\frac{b_{1}}{b_{m}}x+\ldots
            +\frac{b_{m-1}}{b_{m}}x^{m-1}+x^{m}.
\end{array}
\]
By the Division Algorithm, there is a unique $q_{1}(x),$ and a
unique $r_{1}(x)$ in $F[x]$ such that
\[
        f(x)=q_{1}(x)g_{1}(x)+r_{1}(x)\text{ whenever }r_{1}(x) = 0\text{ or }
                \deg r_{1}(x)<\deg g_{1}(x),
\]
that is
\[
\begin{array}{lcl}
    f(x) &=&q_{1}(x)(\beta _{0}+\beta _{1}x+\ldots
            +\beta_{m-1}x^{m-1}+x^{m})+r_{1}(x)\\
        &=&q_{1}(x)(\frac{b_{0}}{b_{m}}+\frac{b_{1}}{b_{m}}x+\ldots
            +\frac{b_{m-1}}{b_{m}}x^{m-1}+x^{m})+r_{1}(x)\\
        &=&q_{1}(x)\frac{1}{b_{m}}(b_{0}+b_{1}x+\ldots
            +b_{m-1}x^{m-1}+b_{m}x^{m})+r_{1}(x)\\
        &=&\left\{\frac{1}{b_{m}}q_{1}(x)\right\} g(x)+r_{1}(x)\\
        &=&q(x)g(x)+r(x).
\end{array}
\]
Thus
\[
q(x)=\frac{1}{b_{m}}q_{1}(x),\text{ and }r(x)=r_{1}(x).
\]
\end{proof}

Although there are many computational algorithms available for obtaining the coefficients of the quotient and remainder polynomials, there is to-date, no explicit algebraic formula that
can be used to compute them. Perhaps, this is because existing algorithms are so computationally efficient. However, computers are incredibly fast nowadays, so that no one is going to complain if an algorithm gives the result in seconds rather than in milliseconds. On the other hand, knowing that a problem can be handled in finite time does not mean we know how long it will take, and so, there will always be a need for efficient algorithms. An interesting discussion of this point, among other things, can be found in \cite{dK02}.

In the present paper, the existence of an algebraic formula will allow the extension of the use of the Division Algorithm into areas that have hitherto not been contemplated. To illustrate this, we consider a novel application where the algebraic formula for the coefficients of the general remainder and quotient polynomial are used to compute a lower Hessenberg-Toeplitz matrices.

\section{The Generalized Quotient Theorem}

If $f(x)=a_{n}x^{n}+a_{n-1}x^{n-1}+\ldots +a_{2}x^{2}+a_{1}x+a_{0}$ where $a_{n}\neq 0$, and
$g(x)=x^{m}+b_{m-1}x^{m-1}+\ldots +b_{2}x^{2}+b_{1}x+b_{0}$, is a monic polynomial with $m<n$ then there exist polynomials $q(x)=d_{n-m}x^{n-m}+d_{n-m-1}x^{n-m-1}+\ldots + d_{2}x^{2}+d_{1}x+d_{0}$ and $r(x)=c_{m-1}x^{m-1}+\ldots +c_{1}x+c_{0}$ such that
\begin{equation}\label{eq-1}
f(x)=g(x)q(x)+r(x).
\end{equation}
Equating the coefficients in (\ref{eq-1}) leads to the following two systems of linear equations. The first system of m linear equations is given by
\begin{equation}\label{eq-2}
c_{k}+\sum_{i+j=k}b_{i}d_{j}=a_{k}
\end{equation}
where $k=0,1,\ldots, m-1$; $i=0,1,\ldots ,m$; and $j=0,1,\ldots ,n-m$. The second system of $n-m$ linear equations is of the form:
\begin{equation}\label{eq-3}
        a_{k}=\sum_{i+j=k}b_{i}d_{j}
\end{equation}
This holds for $k=m,m+1,\ldots ,n$, the ranges for $i$ and $j$ remaining unchanged. Writing (\ref{eq-3}) out, we get
\begin{equation}\label{eq-4}
\begin{array}{lcl}
b_{2m-n}d_{n-m}+b_{2m-n+1}d_{n-m-1}+\ldots +b_{m-1}d_{1}+d_0&=&a_{m},\\
b_{2m-n+1}d_{n-m}+b_{2m-n+2}d_{n-m-1}+\ldots +d_1&=&a_{m+1},\\
&\vdots&\\
b_{m-1}d_{n-m}+d_{n-m-1}&=&a_{n-1},\\
d_{n-m}&=&a_{n},
\end{array}
\end{equation}
with the proviso that $b_{m}=1$ and $b_{i}=0$ when $i$ is
negative, $b_{i}$ will have negative subscript when
$m<\frac{n}{2}$. The system of linear equations (\ref{eq-4}) allows us to compute
$d_{n-m},d_{n-m-1},\ldots ,d_{0}$ recursively, giving,
\begin{equation}\label{d}\small
\begin{array}{lcl}
d_{n-m}&=&a_{n},\\
d_{n-m-1}&=&a_{n-1}-b_{m-1}a_{n},\\
&\vdots&\\
d_{1}&=&a_{m+1}-b_{2m-n+1}a_{n}-b_{2m-n+2}(a_{n-1}-b_{m-1}a_{n})-\ldots -b_{m-1}d_{2},\\
d_{0}&=&a_{m}-b_{2m-n}a_{n}+b_{2m-n+1}(a_{n-1}+b_{m-1}a_{n})+\ldots \\
&& + ~b_{m-1}(a_{+1}-b_{2m-n+1}a_{n}-b_{2m-n+2}(a_{n-1}-b_{m-1}a_{n})- \ldots-\ b_{m-1}d_{2}).
\end{array}
\end{equation}
The coefficients $d_{i}$, $i=0,1,\ldots ,n-m$, can be written in terms of the elements of the linear recurrent sequence $\{s_{n}\}$, where
\[
s_{n+r}=\alpha _{0}s_{n}+\alpha_{1}s_{n+1}+\ldots +\alpha _{r-1}s_{n+r-1},
\]
for positive integers $n$ and $r$, and $\alpha _{i}\in F,$\
$i=0,1,\ldots ,r-1.$\ The result is summarized in Theorem
\ref{th-1} below.

\begin{thm}\label{th-1}
If $f(x)=a_{n}x^{n}+a_{n-1}x^{n-1}+\ldots +a_{2}x^{2}+a_{1}x+a_{0}$ where  $a_{n}\neq 0$, and $g(x)=x^{m}+b_{m-1}x^{m-1}+\ldots +b_{2}x^{2}+b_{1}x+b_{0}$, are polynomials in $F[x]$ with $m<n$, then, the quotient, on dividing $f(x)$\ by $g(x)$ is $q(x)=d_{n-m}x^{n-m}+d_{n-m-1}x^{n-m-1}+\ldots +d_{2}x^{2}+d_{1}x+d_{0}$ where
\begin{equation}\label{eq5}
d_{n-m-k} = \sum_{i=0}^{k}s_{k+1-i}a_{n-i}
\end{equation}
and $\{s_{n}\}$ is the linear recurrent sequence, $s_{1} = 1$ and     $s_{r}=\displaystyle\sum_{i=1}^{r-1}b_{m-i}s_{r-i}$ for $r=2,3,\ldots $.
\end{thm}

\begin{proof} The theorem can be proved by mathematical induction on $k$. It is obviously true for $k=0$.

We assume that it is true for all $d_{n-m-\nu }, ~\nu \le k< n-m$, so that
\begin{equation}\label{eq6}
d_{n-m-\nu }=\sum_{i=0}^{\nu }s_{\nu +1-i}a_{n-i}, \qquad \text{for all}\quad \nu \le k.
\end{equation}
From (\ref{eq-4}), we have
\[
a_{n-\nu -1}=-b_{m}d_{n-m-\nu -1}+b_{m-1}d_{n-m-\nu}+b_{m-2}d_{n-m-\nu +1} +\ldots +b_{m-\nu +1}d_{n-m},
\]
$b_{m}=1$,\ so that
\[
d_{n-m-\nu -1}=a_{n-\nu-1}-b_{m-1}d_{n-m-\nu }-b_{m-2}d_{n-m-\nu +1}+\ldots -b_{m-\nu-1}d_{n-m}.
\]
Substituting from (\ref{eq6}) we get
\[
d_{n-m-\nu -1}=a_{n-\nu -1}-b_{m-1}\sum_{i=0}^{\nu }s_{\nu+1-i}a_{n-i}-b_{m-2}
\sum_{i=0}^{\nu -1}s_{\nu -i}a_{n-i}-\ldots-b_{m-\nu -1}s_{1}a_{n}.
\]
Expanding the summation terms and regrouping we have
\[
d_{n-m-\nu -1}=\sum_{i=0}^{\nu +1}s_{\nu+2-i}a_{n-i}.
\]
Hence (\ref{eq5}) is true for all $k\leq n-m$.
\end{proof}

By Lemma \ref{lem-S1.1}, we can now extend the above result to the case where $b_{m}$ may not equal $1$.

\begin{cor}[Generalized Quotient Theorem]\label{co-1}
If \ $f(x)=a_{n}x^{n}+a_{n-1}x^{n-1}+\ldots +a_{1}x+a_{0}$, $g_{1}(x)=b_{m}x^{m}+b_{m-1}x^{m-1}+\ldots +b_{1}x+b_{0}$ are polynomials in $F[x]$\  of degree $n$ and $m$ respectively where $n>m$, the quotient when $f(x)$ is divided by $g_{1}(x)$ is given by
\[
q_{1}(x)=\bar{d}_{n-m}x^{n-m}+\bar{d}_{n-m-1}x^{n-m-1}+\ldots +\bar{d}_{1}x +\bar{d}_{0}
\]
where
$\bar{d}_{n-m-k}=\frac{1}{b_{m}}\sum_{j=0}^{k}t_{k+1-j}a_{n-j}$\
and $t_{1}=\frac{1}{b_{m}},$
$t_{r}=\frac{1}{b_{m}}\sum_{i=1}^{r-1}b_{m-i}t_{r-i}$ for $r=2,3,\ldots $.
\end{cor}

\begin{proof}
By letting $g(x)=x^{m}+b_{m-1}^{\prime}x^{m-1}+b_{m-2}^{\prime }x^{m-2}+\ldots +b_{1}^{\prime
}x+b_{0}^{\prime }$\ where $b_{m-i}^{\prime }=b_{m-i}/b_{m},$\ the quotient, when $f(x)$ is divided by $g(x),$ is given by
\[
q(x)=d_{n-m}^{\prime}x^{n-m}+d_{n-m-1}^{\prime }x^{n-m-1}+\ldots +d_{1}^{\prime}x+d_{0}^{\prime}
\]
where we note that the sequence $\{t_{k}\}$ is identical to $\{s_{k}\}$ of Theorem \ref{th-1} when $b_{m} = 1$ and hence $d_{n-m-k}^{\prime}=\sum_{i=0}^{k}t_{k+1-i}a_{n-i},$ $t_{1}=\frac{1}{b_{m}},$~ and $t_{r}=\sum_{i=1}^{r-1}b_{m-i}^{\prime }t_{r-i}=\frac{1}{b_{m}}\sum_{i=1}^{r-1}b_{m-i}t_{r-i}.$ Now $g(x)=g_{1}(x)/b_{m}$ so that the quotient $q_{1}(x),$ when $f(x)$\ is divided by $g_{1}(x)$ is $q(x)/b_{m}$ where $\bar{d}_{n-m-k}=\frac{1}{b_{m}}\sum_{j=0}^{k}t_{k+1-j}a_{n-j}.$
\end{proof}

Since the coefficients $a_{0}, a_{1}, \dots, a_{m-1}$ of $f(x)$ do not contribute to the quotient polynomial, the following corollary may be stated without proof.

\begin{cor}\label{col2}
If \ $f(x)=a_{n}x^{n}+a_{n-1}x^{n-1}+\ldots +a_{1}x+a_{0}$, $g(x)= b_{m}x^{m}+b_{m-1}x^{m-1}+\ldots +b_{1}x+b_{0}$ are polynomials in $F[x]$\  of degree $n$ and $m$ respectively where $n>m$, suppose $q(x)$ is the quotient on dividing $f(x)$ by $g(x)$, and
$\hat{q}(x)$ is the quotient on dividing $\hat{f}(x)=a_{n}x^{n}+a_{n-1}x^{n-1}+\ldots +a_{m}x^{m} $ by $g(x)$ then $q(x)=\hat{q}(x)$.
\end{cor}

\section{The Generalized Remainder Theorem} \label{sec4}

The remainder $r(x)$ can now be obtained from (\ref{eq-2}). Since the coefficients $c_{k}$ depends on the coefficients of the quotient polynomial, the most general expression for the remainder can be obtained when the divisor is a general polynomial of degree $m<n$. The result is given by the next theorem.

\begin{thm}[Generalized Remainder Theorem]\label{th2}
If $f(x)=a_{n}x^{n}+a_{n-1}x^{n-1}+\ldots +a_{2}x^{2}+a_{1}x+a_{0}$, $g(x)=b_{m}x^{m}-b_{m-1}x^{m-1}-\ldots -b_{2}x^{2}-b_{1}x-b_{0}$ where $a_{n}\neq 0$ and $b_{m}\neq 0$, are polynomials in $F[x]$ and suppose that $m\leq n$, then the remainder
on dividing $f(x)$ by $g(x)$ is
\begin{equation}\label{GR}
r(x)=\sum_{k=0}^{m-1}\left(a_{k}+\frac{1}{b_{m}}\sum_{i+j=k}b_{i} \sum_{\nu=0}^{n-m-j}t_{n-m-j+1-\nu }a_{n-\nu }\right)x^{k}
\end{equation}
where $\{t_{r}\}$ is the linear recurrent sequence, $t_{1}=\frac{1}{b_{m}}$ and $t_{r}=\frac{1}{b_{m}}\sum_{i=1}^{r-1}b_{m-i}t_{r-i}$ for $r=2,3,\ldots .$
\end{thm}

\begin{proof}
From the system of linear equations (\ref{eq-2}),
\[
c_{k}=a_{k}+\sum_{i+j=k}b_{i}d_{j}
\]
where $k=0,1,\ldots, m-1$; $i=0,1,\ldots ,m$; and $j=0,1,\ldots,n-m$.

Using the result of Corollary \ref{co-1}, and putting $n-m-k=j$ there, we get
\[
d_{j}=\frac{1}{b_{m}}\sum_{\nu=0}^{n-m-j}t_{n-m-j+1-\nu }a_{n-\nu }.
\]
Hence
\[
r(x)=\sum_{k=0}^{m-1}\left(a_{k}+\frac{1}{b_{m}}\sum_{i+j=k}b_{i}
\sum_{\nu=0}^{n-m-j}t_{n-m-j+1-\nu }a_{n-\nu}\right)x^{k}.
    \]
\end{proof}

We note, in passing, that $t_{n-m-j+1-\nu}$ always has positive suffices, by construction.

\section{A Lower Hessenberg-Toeplitz matrix}

Toeplitz matrices is a class of matrices, whose elements along each diagonal are the same constants, see \cite[p.27]{hR85} and \cite[p.1]{xW00}. Elouafi and Aiat Hadj \cite[pp.177-178]{EH09}, defined the lower Hessenberg matrix as follows: The matrix $\mathbf{H} = (h_{ij})_{1\le i,j\le n}$ is called lower Hessenberg matrix if $h_{ij} =0$ for $i<j+1,$ and assume that all elements of the super diagonal of the lower Hessenberg matrix to be non-zero, they given a recursion formula for the characteristic polynomial of the matrix $\mathbf{H}$.
Now we defined the lower Hessenberg-Toeplitz matrices of order $k$, which have the form
\[
\Delta_{k}=\left|\begin{array}{lllll}
    -b_{m-1}    &~~b_{m}    &~~0          &\dots   &~~0\\
    -b_{m-2}    &-b_{m-1}   &~~b_{m}      &\ddots  &~~\vdots\\
    \quad\vdots &\quad\vdots&~~\ddots     &\ddots  &~~0\\
    -b_{m-k+1}  &-b_{m-k+2} &             &\ddots  &~~b_{m}\\
    -b_{m-k}    &-b_{m-k+1} &-b_{m-k+2}   &\dots   &-b_{m-1}
    \end{array}\right|
\]
where $b_{m} \ne 0.$ In this paper we attempt to find the generalized remainder and quotient theorems which related to the determinants of Hessenberg-Toeplitz matrices of this form.

By this second system of linear equations (\ref{eq-4}) where $b_{m} \ne 0$ , we get a matrix equation $\mathbf{H}\mathbf{d}=\mathbf{a}$,

\begin{equation}\label{eq3}
\left(\begin{array}{lllll}
        b_{2m-n}    &b_{2m-n+1} & \ldots    & b_{m-1}  &b_{m}\\
        b_{2m-n+1}  &           &\iddots    & b_{m}    &0\\
        \vdots      &\iddots    &\iddots    &\iddots   &\vdots \\
        b_{m-1}     & b_{m}     &\iddots    &          &\vdots \\
        b_{m}       & 0         & \dots     & \dots    &0 \\
\end{array}
\right)
\left(
\begin{array}{c}
        d_{n-m}\\
        d_{n-m-1}\\
        \vdots\\
        d_{1}\\
        d_{0}\\
\end{array}
\right)
=
\left(
\begin{array}{c}
        a_{m}\\
        a_{m+1}\\
        \vdots\\
        a_{n+1}\\
        a_{n}\\
\end{array}
\right).
\end{equation}
The square backward upper triangular matrix $\mathbf{H}$ of order $k=n-m+1$ is a Hankel matrix with all entries on the backward diagonal as the leading coefficient $b_{m}$ of $g(x)$. So for some $b_{j} = 0$ for $m<\frac{n}{2}$. If then $m<\frac{n}{2}$ for some subscript $j$ of the entries $b_{j}$ in $\mathbf{H}$ is a negative we define $b_{j}=0$ if $j<0$, in this
case the matrices equation (\ref{eq3}) is of the form:
\begin{equation}\label{eq4}
    \left(\begin{array}{ccccccc}
        0           &\dots     &0          & b_{0}    & \dots      &b_{m-1}    &b_{m}\\
        \vdots      &\iddots     &\iddots     &          &\iddots      &\iddots     &0\\
        0           &\iddots     &           &\iddots    &\iddots      &\iddots     &\vdots\\
        b_{0}       &           &\iddots     &\iddots    &\iddots      &           &\vdots\\
        \vdots      & \iddots    &\iddots     &\iddots    &            &           &\vdots\\
        b_{m-1}     & b_{m}     &\iddots     &          &            &           &\vdots \\
        b_{m}       & 0         &\dots      &\dots     &\dots       &\dots      &0\\
    \end{array}\right)
    \left(\begin{array}{c}
        d_{n-m}\\
        d_{n-m-1}\\
        \vdots\\
        \vdots\\
        \vdots\\
        d_{1}\\
        d_{0}\\
    \end{array}\right)
    =
    \left(\begin{array}{c}
        a_{m}\\
        a_{m+1}\\
        \vdots\\
        \vdots\\
        \vdots\\
        a_{n-1}\\
        a_{n}\\
    \end{array}\right).
\end{equation}

Now, let $\mathbf{W}$ be a bordered matrix obtained from $\mathbf{H}$ and $\mathbf{a}$ see \cite[p.417]{aJ77}, defined as follows:
\begin{equation}\label{eq7}
    \mathbf{W} =\left(\begin{array}{cc}
        \mathbf{H}       & \mathbf{a}\\
        \mathbf{x}^{T}   & 0
        \end{array}\right)
\end{equation}
where the matrix $\mathbf{H}$ and $\mathbf{a}$ is define in (\ref{eq3}) and
\[
    \mathbf{x}^{T}=\left(\begin{array}{cccccc}
            x^{n-m}  &x^{n-m}   &\dots   &x^{2} &x   &1
        \end{array}\right).
\]
Thus $\mathbf{W}$ is a square matrix of order $t=k+1=n-m+2$.

    From (\ref{eq7})
\begin{equation}\label{eq8}
    \det(\mathbf{W})=\left|\begin{array}{cc}
        \mathbf{H}       & \mathbf{a}\\
        \mathbf{x}^{T}   & 0\\
        \end{array}\right|
        =
        \left|\begin{array}{llllll}
        b_{2m-n}    &b_{2m-n+1} & \ldots    & b_{m-1}  &b_{m}   &a_{m}\\
        b_{2m-n+1}  &           & \iddots   & b_{m}    &0       &a_{m+1}\\
        \vdots      & \qq\iddots   & \iddots   & \iddots  &\vdots  &\vdots \\
        b_{m-1}     & b_{m}     & \iddots   &          &\vdots  &a_{n-1}\\
        b_{m}       & 0         & \ldots    & \ldots   &0       &a_{n}\\
        x^{n-m}     & x^{n-m-1} & \ldots    & x        &1       &0
        \end{array}\right|.
\end{equation}
Interchange the last column with the consecutive columns of the
matrix $\mathbf{W}$ until to the first column is reached and similarly
interchange the last row of the matrix $\mathbf{W}$ with the consecutive
row until to the top row is reached. The sum of the number of
interchange is even. Thus the new matrix, namely $\mathbf{T}$ form from the matrix
$\mathbf{W}$ has order $t$ and has the same determinant as $\mathbf{W}$. That is
\begin{equation}\label{eq-T}
\mathbf{T}=\left(\begin{array}{llllll}
            0       & x^{n-m}   & x^{n-m-1}     & \dots     &x          &1\\
            a_{m}   & b_{2m-n}  & b_{2m-n+1}    & \dots     &b_{m-1}    &b_{m}\\
            a_{m-1} & b_{2m-n+1}&               & \iddots   &b_{m}      &0 \\
            \vdots  & \vdots    &\qq \iddots       & \iddots   &\iddots    &\vdots \\
            a_{n-1} & b_{m-1}   & b_{m}         & \iddots   &           &\vdots \\
            a_{n}   & b_{m}     & 0             & \dots     &\dots      &0
        \end{array}\right),
\end{equation}
and
\begin{equation}\label{eq9}
\det(\mathbf{W})
=\left|\begin{array}{llllll}
0       & x^{n-m}   & x^{n-m-1}     & \dots     &x          &1\\
a_{m}   & b_{2m-n}  & b_{2m-n+1}    & \dots     &b_{m-1}    &b_{m}\\
a_{m+1} & b_{2m-n+1}&               & \iddots   &b_{m}      &0 \\
\vdots  & \vdots    &\qq\iddots       & \iddots   &\iddots    &\vdots \\
a_{n+1} & b_{m-1}   & b_{m}         & \iddots   &           &\vdots \\
a_{n}   & b_{m}     & 0             & \dots     &\dots      &0
\end{array}\right|
=\det(\mathbf{T}).
\end{equation}

    Now consider a backward identity matrix (or anti-identity matrix) of
order $t=n-m+2$, $\mathbf{P}_{t}$ say,

\[
    \mathbf{P}_{t}=\mathbf{P}_{n-m+2}=\left(\begin{array}{llllll}
    0       &\dots    &\dots      &\dots     &0        &1\\
    \vdots  &         &           &\iddots    &1        &0\\
    \vdots  &         & \iddots    &\iddots    &\iddots   &\vdots\\
    \vdots  &\iddots   & \iddots    &\iddots    &         &\vdots\\
    0       &1        & \iddots    &          &         &\vdots\\
    1       &0        & \dots     &\dots     &\dots    &0\\
    \end{array}\right)_{(n-m+2,n-m+2)}.
\]

The following lemma and its proof are well known. We present
this material here because we refer to it a few times in the next
two theorems.

\begin{lem}\label{lem-1.1}
If $\mathbf{P}_{t}$ (as above) is the back word identity matrix of order $t$, then
\begin{equation}\label{eq91}
        \det(\mathbf{P}_{t})=(-1)^{\frac{t(t-1)}{2}}.
\end{equation}
\end{lem}

Now consider the product of the matrix $\mathbf{T}$ from \eqref{eq-T} and $\mathbf{P}_{t}$ defined as above, \cite{hR85}, assert that $\mathbf{P}_{t}\mathbf{T}=\mathbf{L}$ (and $\mathbf{TP}_{t}=\mathbf{U}$) where
$\mathbf{L}$ is a lower Hessenberg-Toeplitz matrix (and $\mathbf{U}$ is an upper
Hessenberg-Toeplitz matrix). In general
\[
\begin{array}{lcl}
    \mathbf{P}_{t}\mathbf{T} &=&\left(\begin{array}{llllll}
    0       &\dots    &\dots      &\dots     &0        &1\\
    \vdots  &         &           &\iddots    &1        &0\\
    \vdots  &         & \iddots    &\iddots    &\iddots   &\vdots\\
    \vdots  &\iddots   & \iddots    &\iddots    &         &\vdots\\
    0       &1        & \iddots    &          &         &\vdots\\
    1       &0        & \dots     &\dots     &\dots    &0\\
    \end{array}\right)\\

    &&\times
    \left(\begin{array}{llllll}
            0       & x^{n-m}       & x^{n-m-1}     & \dots     &x          &1\\
            a_{m}   & b_{2m-n}      & b_{2m-n+1}    & \dots     &b_{m-1}    &b_{m}\\
            a_{m-1} & b_{2m-n+1}    &               & \iddots    &b_{m}      &0 \\
            \vdots  & \vdots        & \qq\iddots        & \iddots    &\iddots     &\vdots \\
            a_{n-1} & b_{m-1}       & b_{m}         & \iddots    &           &\vdots \\
            a_{n}   & b_{m}         & 0             & \dots     &\dots      &0
        \end{array}\right)\\
    &=&\left(\begin{array}{llllll}
            a_{n}   & b_{m}         & 0             & \dots     &\dots     &0\\
            a_{n-1} & b_{m-1}       & b_{m}         & \ddots    &          &\vdots \\
            \vdots  & \vdots        & \ddots        & \ddots    &\ddots    &\vdots \\
            a_{m-1} & b_{2m-n+1}    &               & \ddots    &b_{m}     &0 \\
            a_{m}   & b_{2m-n}      & b_{2m-n+1}    & \dots     &b_{m-1}   &b_{m}\\
            0       & x^{n-m}       & x^{n-m-1}     & \dots     &x         &1
        \end{array}\right)=: \mathbf{H}_{t}.
\end{array}\]
Let $\Delta_{i}:=\det(\mathbf{H}_{i})$. Expanding the determinant of
Hessenberg-Toeplitz matrix, $\det(\mathbf{H}_{t})$
\begin{equation}\label{eq10}
    \det(\mathbf{H}_{t})=\left|\begin{array}{llllll}
            a_{n}   & b_{m}         & 0             & \dots     &\dots      &0\\
            a_{n-1} & b_{m-1}       & b_{m}         & \ddots    &           &\vdots \\
            \vdots  & \vdots        & \ddots        & \ddots    &\ddots     &\vdots \\
            a_{m-1} & b_{2m-n+1}    &               & \ddots    &b_{m}      &0 \\
            a_{m}   & b_{2m-n}      & b_{2m-n+1}    & \dots     &b_{m-1}    &b_{m}\\
            0       & x^{n-m}       & x^{n-m-1}     & \dots     &x          &1
        \end{array}\right|,
\end{equation}
by the last row. Thus
\begin{equation}\label{eq11}
\begin{array}{lcl}
\det(\mathbf{H}_{t}) &=&1\Delta _{t-1}-xb_{m}\Delta_{t-2}+x^{2}b_{m}^{2}\Delta _{t-3}
-x^{3}b_{m}^{3}\Delta_{t-4}+\ldots-x^{n-m}b_{m}^{n-m}\Delta_{1}+0\\
&=&\displaystyle\sum_{i=2}^{t}(-1)^{t-i}x^{t-i}b_{m}^{t-i}\Delta_{i-1}
\end{array}
\end{equation}
where
\[
\Delta_{1}=\left|a_{n}\right|,~~
\Delta_{2}=\left|\begin{array}{ll}
            a_{n}       &b_{m}\\
            a_{n-1}     &b_{m-1}
            \end{array}\right|,~~
\Delta_{3}=\left|\begin{array}{lll}
            a_{n}     &b_{m}    &0\\
            a_{n-1}   &b_{m-1}  &b_{m}\\
            a_{n-2}   &b_{m-2}  &b_{m-1}
            \end{array}\right|,
\dots
\]
\begin{equation}
\Delta_{t-1}=\left|\begin{array}{lllll}
    a_{n}      &b_{m}       &0            &\dots   &0\\
    a_{n-1}    &b_{m-1}     &b_{m}        &\ddots  &\vdots\\
    \vdots     &\vdots      &\ddots       &\ddots  &0\\
    a_{m+1}    &b_{2m-n+1}  &             &\ddots  &b_{m}\\
    a_{m}      &b_{2m-n}    &b_{2m-n+1}   &\dots   &b_{m-1}
    \end{array}\right|.\label{eq-Delta t-1}
\end{equation}

\medskip

\begin{thm}\label{thm-1.2}
If \;$\mathbf{T}$ and $\mathbf{H}_{t}$ are matrices of order $t$ defined as above then
\begin{equation}\label{eq111}
        \det(\mathbf{T}) =\det(\mathbf{W}) = (-1)^{\frac{t(t-1)}{2}}\det(\mathbf{H}_{t})
\end{equation}
\end{thm}

\begin{proof}
Since $\det(\mathbf{H}_{t})=\det (\mathbf{P}_{t}\mathbf{T})=\det(\mathbf{P}_{t})\det(\mathbf{T})$ implies that $\det(\mathbf{T})=\displaystyle\frac{\det(\mathbf{H}_{t})}{\det(\mathbf{P}_{t})}$. From (\ref{eq9}) $\det(\mathbf{T})=\det(\mathbf{W})$ we have
$\det(\mathbf{W})=\displaystyle\frac{\det(\mathbf{H}_{t})}{\det(\mathbf{P}_{t})}.$ Since, $\det(\mathbf{P}_{t})=(-1)^{\frac{t(t-1)}{2}}$ , by Lemma \ref{lem-1.1}. Thus \eqref{eq111}
true.
\end{proof}

\begin{cor}\label{cor-1.3}
$\det(\mathbf{W}) = (-1)^{\frac{t(t-1)}{2}}\sum_{i=2}^{t}(-1)^{t-i}x^{t-i}b_{m}^{t-i}\Delta _{i-1}$.
\end{cor}

\begin{proof}
By \eqref{eq111} and \eqref{eq11}.
\end{proof}

\begin{thm}\label{thm-1.4}
The determinant of the Hankel matrix,
\[
    \mathbf{H} = \left(\begin{array}{lllll}
        b_{2m-n}        &b_{2m-n+1} &\dots   &b_{m-1}    &b_{m}\\
        b_{2m-n+1}      &           &\iddots &b_{m}      &0\\
        \vdots          &\qq\iddots    &\iddots &\iddots    &\vdots\\
        b_{m-1}         &b_{m}      &\iddots &           &\vdots\\
        b_{m}           &0          &\dots   &\dots      &0\\
        \end{array}\right)_{(t-1,t-1)}.
\]
is $\det(\mathbf{H}) =(-1)^{\frac{(t-1)(t-2)}{2}}b_{b}^{t-1}$.
\end{thm}

\begin{proof}
Since
\[\small
\begin{array}{lcl}
\mathbf{H}\mathbf{P}_{t-1} &=& \left(\begin{array}{lllll}
        b_{2m-n}        &b_{2m-n+1} &\dots   &b_{m-1}    &b_{m}\\
        b_{2m-n+1}      &           &\iddots &b_{m}      &0\\
        \vdots          &\qq\iddots    &\iddots &\iddots    &\vdots\\
        b_{m-1}         &b_{m}      &\iddots &           &\vdots\\
        b_{m}           &0          &\dots   &\dots      &0\\
        \end{array}\right)
    \left(\begin{array}{lllll}
        0       &\dots  &\dots  &0          &1\\
        \vdots  &       &\iddots &1          &0\\
        \vdots  &\iddots &\iddots &\iddots     &\vdots\\
        0       &1      &\iddots  &          &\vdots\\
        1       &0      &\dots  &\dots      &0\\
        \end{array}\right)\smallskip\\
&=&
    \left(\begin{array}{lllll}
        b_{m}   &b_{m-1}    &\dots          &b_{2m-n+1} &b_{2m-n}\\
        0       &b_{m}      &\ddots         &           &b_{2m-n+1}\\
        \vdots  &\ddots     &\ddots         &\ddots     &\vdots\\
        \vdots  &           &\ddots         &b_{m}      &b_{m-1}\\
        0   &\dots      &\dots          &0          &b_{m}\\
        \end{array}\right)=:\mathbf{U}.
\end{array}
\]
Since $\det (\mathbf{H}\mathbf{P}_{t-1})=\det(\mathbf{H})\det(\mathbf{P}_{t-1})=\det(\mathbf{U}) =b_{m}^{t-1}.$ Thus,
\[
\det(\mathbf{H})=\frac{\det \mathbf{U}}{\det(\mathbf{P}_{t-1})}.
\]
By Lemma \ref{lem-1.1},
\[
\det(\mathbf{P}_{t-1}) =(-1)^{\frac{(t-1)[(t-1)-1]}{2}}=(-1)^{\frac{(t-1)(t-2)}{2}}.
\]
Therefore:
\[
    \det(\mathbf{H}) =\frac{\det\mathbf{U}}{\det(\mathbf{P}_{t-1})}=\frac{b_{m}^{t-1}}{(-1)^{\frac{(t-1)(t-2)}{2}}}
    =(-1)^{\frac{(t-1)(t-2)}{2}}b_{m}^{t-1}.
\]
The theorem was proved.
\end{proof}

\begin{thm}\label{thm-1.5}
The quotient of the polynomial $f(x)$ by ~$g(x)$ is $q(x)=-\displaystyle\frac{\det(\mathbf{W})}{\det(\mathbf{H})}.$
\end{thm}

\begin{proof}
From (\ref{eq8})
\[\begin{array}{lcl}
    \det(\mathbf{W}) &=& \left|\begin{array}{ll}
            \mathbf{H}       &\mathbf{a}\\
            \mathbf{x}^{T}   &0\\
            \end{array}\right|\smallskip\\
            &=&
    \left|\begin{array}{llllll}
        b_{2m-n}   &b_{2m-n+1} &\dots   &b_{m-1}     &b_{m}      &a_{m}\\
        b_{2m-n+1} &           &\iddots  &b_{m}       &0          &a_{m+1}\\
        \vdots     &\qq\iddots     &\iddots  &\iddots      &\vdots     &\vdots\\
        b_{m-1}    &b_{m}      &\iddots  &            &\vdots     &a_{n-1}\\
        b_{m}      &0          &\dots   &\dots       &0          &a_{n}\\
        x^{n-m}    &x^{n-m-1}  &\dots   &x           &1          &0\\
        \end{array}\right|
\end{array}
\]
and J.W. Archbold \cite[p.417]{aJ77}, assert that:
\begin{equation}\label{eq12}
\det(\mathbf{W})
=-\mathbf{x}^{T}\text{adj}(\mathbf{H})\mathbf{a}+0\det(\mathbf{H})\\
=-\mathbf{x}^{T}\text{adj}(\mathbf{H})\mathbf{a}.
\end{equation}
Consider \eqref{eq3}, the matrix equation $\mathbf{H}\mathbf{d}=\mathbf{a}$, we get
$\mathbf{d}=\mathbf{H}^{-1}\mathbf{a}$, and Theorem \ref{thm-1.4}, show that $\mathbf{H}$ is
nonsingular matrix, thus $\mathbf{H}^{-1}$ exists. From the inverse formula
            $\mathbf{H}^{-1}=\displaystyle\frac{1}{\det(\mathbf{H})}\text{adj}(\mathbf{H})$
we get:
\[
\mathbf{d}=\frac{1}{\det(\mathbf{H})}\text{adj}(\mathbf{H})\mathbf{a}.
\]
That is:
\[
\mathbf{x}^{T}\mathbf{d}=\mathbf{x}^{T}(\frac{1}{\det(\mathbf{H})}\text{adj}(\mathbf{H}) \mathbf{a})
             =\frac{1}{\det(\mathbf{H})}(\mathbf{x}^{T}\text{adj}(\mathbf{H})\mathbf{a})
\]
Since $q(x)=\mathbf{x}^{T}\mathbf{d}$ so that:
\begin{equation}\label{eq13}
q(x)=\frac{1}{\det(\mathbf{H})}(\mathbf{x}^{T}\text{adj}(\mathbf{H})\mathbf{a})
\end{equation}
From \eqref{eq12}
\[
\begin{array}{lcl}
\det(\mathbf{W}) &=&-\mathbf{x}^{T}\text{adj}(\mathbf{H})\mathbf{a}\\
&=&-\det(\mathbf{H})\left(\displaystyle\frac{1}{\det(\mathbf{H})}
       (\mathbf{x}^{T}\text{adj}(\mathbf{H}) \mathbf{a})\right)\\
&=&-\det(\mathbf{H})q(x).
\end{array}
\]
Therefore $q(x)=-\displaystyle\frac{\det(\mathbf{W})}{\det(\mathbf{H})}$.
\end{proof}

\begin{thm}\label{thm-2.1}
If $f(x)=a_{n}x^{n}+a_{n-1}x^{n-1}+\ldots +a_{2}x^{2}+a_{1}x+a_{0}$ and $g(x)=b_{m}x^{m}-b_{m-1}x^{m-1}-\ldots -b_{2}x^{2}-b_{1}x-b_{0}$  where $a_{n}\neq 0$, and suppose that $m=\deg g(x)\leq \deg f(x)=n$ then the quotient on dividing $f(x)$ by $g(x)$ is
\[
    q(x)=\sum_{j=0}^{n-m}\left( (-1)^{n-m-j}b_{m}^{j-(n-m+1)}
                                \Delta_{(n-m+1)-j}\right) x^{j},
\]
where $\Delta_{(n-m+1)-j}$ define as in \eqref{eq-Delta t-1}.
\end{thm}

\begin{proof}
From Corollary \ref{cor-1.3}, $\det(\mathbf{W})=(-1)^{\frac{t(t-1)}{2}}\sum_{i=2}^{t}(-1)^{t-i}x^{t-i}b_{m}^{t-i}\Delta
_{i-1}$  and Theorem \ref{thm-1.5}, $q(x)=-\frac{\det(\mathbf{W})}{\det(\mathbf{H})}$, and
Theorem \ref{thm-1.4}, $\det(\mathbf{H})=(-1)^{\frac{(r-1)(t-2)}{2}}b_{m}^{t-1}$ we
have:
\[
\begin{array}{lcl}
    q(x)&=&-\displaystyle\frac{(-1)^{\frac{t(t-1)}{2}}\sum_{i=2}^{t}(-1)^{t-i}
            x^{t-i}b_{m}^{t-i}\Delta_{i-1}}{(-1)^{\frac{(r-1)(t-2)}{2}}b_{m}^{t-1}}\\
        &=&(-1)^{[1+\frac{t(t-1)}{2}-\frac{(r-1)(t-2)}{2}]}b_{m}^{-(t-1)}
            \displaystyle\sum_{i=2}^{t}(-1)^{t-i}x^{t-i}b_{m}^{t-i}\Delta_{i-1}\\
        &=&(-1)^{t}\displaystyle\sum_{i=2}^{t}(-1)^{t-i}b_{m}^{t-i-(t-1)}\Delta_{i-1}x^{t-i}\\
        &=&\displaystyle\sum_{i=2}^{t}(-1)^{t}(-1)^{t-i}b_{m}^{1-i}\Delta_{i-1}x^{t-i}\\
        &=&\displaystyle\sum_{i=2}^{t}(-1)^{i}b_{m}^{1-i}\Delta_{i-1}x^{t-i}.\\
\end{array}
\]
Therefore $q(x)=\sum_{i=2}^{t}(-1)^{i}b_{m}^{1-i}\Delta
_{i-1}x^{t-i}=d_{0}+d_{1}x+d_{2}x^{2}+\ldots +d_{n-m}x^{n-m}$, so
that,
\[\begin{array}{lcl}
    d_{0} &=&(-1)^{t}b_{m}^{1-t}\Delta _{t-1}\\
            &=&(-1)^{(n-m+2)}b_{m}^{1-(n-m+2)}\Delta_{(n-m+2)-1}\\
            &=&(-1)^{(n-m+2)}b_{m}^{m-n-1}\Delta_{n-m+1},\\
    d_{1} &=&(-1)^{t-1}b_{m}^{1-(t-1)}\Delta_{(t-1)-1}\\
            &=&(-1)^{t-1}b_{m}^{2-t}\Delta_{t-2},\\
    d_{2}&=&(-1)^{t-2}b_{m}^{1-(t-2)}\Delta_{(t-2)-1}\\
            &=&(-1)^{t-2}b_{m}^{3-t}\Delta_{t-3},\\
            &\vdots& \\
    d_{n-m-1}&=&(-1)^{t-(n-m-1)}b_{m}^{1-(t-(n-m-1))}\Delta_{(t-(n-m-1))-1}\\
            &=&(-1)^{(n-m+2)-(n-m-1)}b_{m}^{1-((n-m+2)-(n-m-1))}\Delta_{((n-m+2)-(n-m-1))-1}\\
            &=&(-1)^{3}b_{m}^{-2}\Delta _{2}\\
     d_{n-m}&=&(-1)^{t-(n-m)}b_{m}^{1-(t-(n-m))}\Delta_{(t-(n-m))-1}\\
            &=&(-1)^{(n-m+2)-(n-m)}b_{m}^{1-((n-m+2)-(n-m))}\Delta_{((n-m+2)-(n-m))-1}\\
            &=&(-1)^{2}b_{m}^{-1}\Delta _{1}.
\end{array}
\]
That is
\begin{equation}\label{eq14}
    d_{j}=(-1)^{(t-j)}b_{m}^{(j+1)-t}\Delta _{t-(j+1)},~~ j=0,1,\dots,n-m.
\end{equation}
Since $q(x)=\sum_{j=0}^{n-m}d_{j}x^{j}$ and $t=n-m+2$ we arrive
at:
\[
\begin{array}{lcl}
    q(x) &=&\displaystyle\sum_{j=0}^{n-m}\left(
            (-1)^{t+j}b_{m}^{j+1-t}\Delta_{t-(j+1)}\right)x^{j}\\
         &=&\displaystyle\sum_{j=0}^{n-m}\left(
            (-1)^{(n-m+2)+j}b_{m}^{j+1-(n-m+2)}\Delta
            _{(n-m+2)-(j+1)}\right)x^{j}.
\end{array}
\]
Therefore: $q(x)=\sum_{j=0}^{n-m}\left(
            (-1)^{n-m+j}b_{m}^{j-(n-m+1)}\Delta_{(n-m+1)-j}\right) x^{j}$.
\end{proof}

\section{Determinant of the Hessenberg-Toeplitz Matrix}

Let $f(x)=a_{n}x^{n}+a_{n-1}x^{n-1}+\ldots +a_{2}x^{2}+a_{1}x+a_{0}$ and $g(x)=b_{m}x^{m}-b_{m-1}x^{m-1}-\ldots -b_{2}x^{2}-b_{1}x-b_{0}$ be polynomials in $F[x]$  degree $n$ and $m$ respectively. In special case, if $b_{m-1}\neq 0$ in $F$, and
\begin{eqnarray}\label{eq17}
    &&a_{n}=-b_{m-1},~a_{n-1}=-b_{m-2},~\dots,~a_{n-m+2}=-b_{1},~a_{n-m+1}=-b_{0},\nonumber\\
    &&a_{n-m}=0,~a_{n-m-1}=0,~\dots,~a_{1}=0,~a_{0}=0.
\end{eqnarray}
then
\[\begin{array}{lcl}
    f(x)&=&0+0x+\ldots
            +0x^{n-m}-b_{m-m}x^{n-m+1}-\dots-b_{m-2}x^{n-1}-b_{m-1}x^{n}\\
        &=&-b_{m-m}x^{n-m+1}-\dots-b_{m-2}x^{n-1}-b_{m-1}x^{n}.
\end{array}
\]
Consider the system of linear equation \eqref{eq-2} is change to
\begin{equation}\label{eq18}
\begin{array}{lcl}
    -b_{2m-n}d_{n-m}-b_{2m-n+1}d_{n-m-1}-\ldots-b_{m-1}d_{1}+b_{m}d_{0}&=&-b_{2m-n-1},\\
    -b_{2m-n+1}d_{n-m}-b_{2m-n+2}d_{n-m-1}-\ldots+b_{m}d_{1}                 &=&-b_{2m-n},\\
                                                                       &\vdots&\\
    -b_{m-1}d_{n-m}+b_{m}d_{n-m-1}                                     &=&-b_{m-2},\\
    b_{m}d_{n-m}                                                       &=&-b_{m-1},
\end{array}
\end{equation}
we define $b_{j}=0$, for $j<0$. Therefore the Hessenberg matrix of
(\ref{eq10}) is change to:
\begin{equation}\label{eq19}
    \mathbf{H}_{t}=\left(\begin{array}{cccccc}
            -b_{m-1}     & ~~b_{m}       & 0             & \dots     &\dots      &0\\
            -b_{m-2}     & -b_{m-1}      & ~~b_{m}       & \ddots    &           &\vdots \\
            \quad\vdots  &\quad \vdots   & \ddots        & \ddots    &\ddots     &\vdots \\
            -b_{2m-n}    & -b_{2m-n+1}   &               & \ddots    &~~b_{m}    &0 \\
            -b_{2m-n-1}  & -b_{2m-n}     & -b_{2m-n+1}   & \dots     &-b_{m-1}   &b_{m}\\
             0           & x^{n-m}       & x^{n-m-1}     & \dots     &x          &1
        \end{array}\right).
\end{equation}
As in \eqref{eq11}, the determinant of this matrix is:
\[
\begin{array}{lcl}
    \det(\mathbf{H}_{t}) &=&1\Delta _{t-1}-xb_{m}\Delta_{t-2}+x^{2}b_{m}^{2}\Delta _{t-3}+\ldots
                    -x^{n-m}b_{m}^{n-m}\Delta _{1}+0\\
                &=&\sum_{i=2}^{t}(-1)^{t-i}x^{t-i}b_{m}^{t-i}\Delta_{i-1}.\\
\end{array}
\]
where
\[
\Delta _{1}=\left|-b_{m-1}\right|, \Delta
_{2}=\left|\begin{array}{ll}
            -b_{m-1}   &~~b_{m}\\
            -b_{m-2}   &-b_{m-1}
            \end{array}\right|,
\Delta _{3}=\left|\begin{array}{lll}
            -b_{m-1}    &~~b_{m}   &~~0\\
            -b_{m-2}    &-b_{m-1}  &~~b_{m}\\
            -b_{m-3}    &-b_{m-2}  &-b_{m-1}
            \end{array}\right|, \dots
\]
\begin{equation}\label{eq20}
\Delta _{t-1}=\left|\begin{array}{lllll}
    -b_{m-1}    &~~b_{m}        &~~0          &\dots   &~~0\\
    -b_{m-2}    &-b_{m-1}       &~~b_{m}      &\ddots  &~~\vdots\\
    \quad\vdots &\quad\vdots    &~~\ddots     &\ddots  &~~0\\
    -b_{2m-n}   &-b_{2m-n+1}    &             &\ddots  &~~b_{m}\\
    -b_{2m-n-1} &-b_{2m-n}      &-b_{2m-n+1}  &\dots   &-b_{m-1}
    \end{array}\right|.
\end{equation}

The following is the main result:

\begin{thm}\label{thm-3.1}
If $\Delta_{k}=\left|\begin{array}{lllll}
    -b_{m-1}    &~~b_{m}    &~~0          &\dots   &~~0\\
    -b_{m-2}    &-b_{m-1}   &~~b_{m}      &\ddots  &~~\vdots\\
    \quad\vdots &\quad\vdots&~~\ddots     &\ddots  &~~0\\
    -b_{m-k+1}  &-b_{m-k+2} &             &\ddots  &~~b_{m}\\
    -b_{m-k}    &-b_{m-k+1} &-b_{m-k+2}   &\dots   &-b_{m-1}
    \end{array}\right|$ then
\[
\Delta_{k}=(-1)^{k}b_{m}^{k}\displaystyle\sum_{i=1}^{k}t_{i}b_{m-(k+1)+i}
=(-1)^{k}b_{m}^{k}\displaystyle\sum_{j=1}^{k-1}t_{k-j}b_{(m-1)-j}
\]
where $\{t_{r}\}$ is the linear recurrent sequence,
$t_{1}=\frac{1}{b_{m}}$ and
$t_{r}=\frac{1}{b_{m}}\sum_{i=1}^{r-1}b_{m-i}t_{r-i}$\ for
$r=2,3,\ldots .$
\end{thm}

\begin{proof}
From \eqref{eq20}
\[
\Delta _{t-1}=\left|\begin{array}{lllll}
    -b_{m-1}    &~~b_{m}        &~~0            &\dots   &~~0\\
    -b_{m-2}    &-b_{m-1}       &~~b_{m}        &\ddots  &~~\vdots\\
    \quad\vdots &\quad\vdots    &~~\ddots       &\ddots  &~~0\\
    -b_{2m-n}   &-b_{2m-n+1}    &               &\ddots  &~~b_{m}\\
    -b_{2m-n-1} &-b_{2m-n}      &-b_{2m-n+1}    &\dots   &-b_{m-1}
    \end{array}\right|,
\]
where $t=n-m+2$. Let $k=n-m+1$. Thus this determinant change to:
\[
\Delta _{k}=\left|\begin{array}{lllll}
    -b_{m-1}    &~~b_{m}    &~~0          &\dots   &~~0\\
    -b_{m-2}    &-b_{m-1}   &~~b_{m}      &\ddots  &~~\vdots\\
    \quad\vdots &\quad\vdots&~~\ddots     &\ddots  &~~0\\
    -b_{m-k+1}  &-b_{m-k+2} &             &\ddots  &~~b_{m}\\
    -b_{m-k}    &-b_{m-k+1} &-b_{m-k+2}   &\dots   &-b_{m-1}
    \end{array}\right|.
\]
According to Theorem \ref{thm-2.1},
\begin{equation}\label{eq21}
    q(x)=\sum_{i=0}^{n-m}\left( (-1)^{n-m-i}b_{m}^{i-(n-m+1)}
            \Delta_{(n-m+1)-i}\right) x^{i}.
\end{equation}
and Corollary \ref{co-1} asserted that
\[
q(x)=\sum_{i=0}^{n-m}\left(\sum_{j=0}^{n-m-i}t_{n-m-i+1-j}a_{n-j}\right) x^{i}.
\]
where $\{t_{r}\}$ is the linear recurrent sequence, $t_{1}=\frac{1}{b_{m}},$ and $t_{r}=\frac{1}{b_{m}}\sum_{i=1}^{r-1}b_{m-i}t_{r-i}$ for $r=2,3,\ldots .$
In this case, by \eqref{eq17}, we get
\begin{equation}\label{eq22}
    q(x)=-\sum_{i=0}^{n-m}\left(\sum_{j=0}^{n-m-i}t_{n-m-i+1-j}b_{(m-1)-j}\right)x^{i}.
\end{equation}

\noindent By the uniqueness of the quotient $q(x) = d_{0} + d_{1}x + \dots + d_{n-m}x^{n-m}$, from \eqref{eq21} and \eqref{eq22} we
have
\[\begin{array}{lcl}
    q(x)&=&\displaystyle\sum_{i=0}^{n-m}\left( (-1)^{n-m-i}b_{m}^{i-(n-m+1)}
            \Delta_{(n-m+1)-i}\right)x^{i}\\
        &=&-\displaystyle\sum_{i=0}^{n-m}\left(\sum_{j=0}^{n-m-i}
        t_{n-m-i+1-j}b_{(m-1)-j}\right)x^{i}.
\end{array}
\]
Equating the coefficients of the quotient polynomial, we get
\begin{eqnarray}\label{eq221}
    d_{0} &=&(-1)^{n-m-0}b_{m}^{0-(n-m+1)}\Delta_{(n-m+1)-0}\nonumber\\
          &=&-\sum_{j=0}^{n-m-0}t_{n-m-0+1-j}b_{(m-1)-j}\\
    d_{1} &=&(-1)^{n-m-1}b_{m}^{1-(n-m+1)}\Delta_{(n-m+1)-1}\nonumber\\
          &=&-\sum_{j=0}^{n-m-1}t_{n-m-1+1-j}b_{(m-1)-j}\nonumber\\
          &\vdots& \nonumber\\
    d_{n-m-1}&=&(-1)^{n-m-(n-m-1)}b_{m}^{(n-m-1)-(n-m+1)}\Delta_{(n-m+1)-(n-m-1)}\nonumber\\
             &=&-\sum_{j=0}^{n-m-(n-m-1)}t_{n-m-(n-m-1)+1-j}b_{(m-1)-j}\nonumber\\
    d_{n-m}  &=&(-1)^{n-m-(n-m)}b_{m}^{(n-m)-(n-m+1)}\Delta_{(n-m+1)-(n-m)}\nonumber\\
             &=&-\sum_{j=0}^{n-m-(n-m-1)}t_{n-m-(n-m)+1-j}b_{(m-1)-j}.\nonumber
\end{eqnarray}
We claim that
\begin{equation}\label{eq23}
\Delta_{k}=(-1)^{k}b_{m}^{k}\sum_{i=1}^{k}t_{i}b_{m-(k+1)+i}.
\end{equation}
Equation \eqref{eq23} can proof by mathematical induction on $k$.

For $k=1$, consider
\[
    d_{n-m}=(-1)^{0}b_{m}^{-1}\Delta_{1}=-\sum_{j=0}^{0}t_{1-1}b_{(m-1)-j}
\]
Since
$(-1)^{0}b_{m}^{-1}\Delta_{1}=\frac{1}{b_{m}}|-b_{m-1}|=-t_{1}b_{m-1}$
implies that
\[
\Delta_{1}=-b_{m}t_{1}b_{m-1}=(-1)^{1}b_{m}^{1}\sum_{i=1}^{1}t_{i}b_{m-(1+1)-i}
\]
Thus the equation \eqref{eq23} is true for $k=1.$

Assume that the equation \eqref{eq23} is true for all order of the determinant in \eqref{eq20} $k<n-m+1$, we must show that it true for $k=n-m+1$. From \eqref{eq221}
\[
d_{0}=(-1)^{n-m}b_{m}^{-(n-m+1)}\Delta_{(n-m+1)}=-\sum_{j=0}^{n-m}t_{n-m+1-j}b_{(m-1)-j}
\]
implies
\[
\begin{array}{lcl}
\Delta_{(n-m+1)}&=&(-1)^{1-n+m}b_{m}^{n-m+1}\sum_{j=0}^{n-m}t_{n-m+1-j}b_{(m-1)-j}\smallskip\\
&=&(-1)^{n-m+1}b_{m}^{n-m+1}\sum_{j=0}^{n-m}t_{n-m+1-j}b_{(m-1)-j}
\end{array}
\]
If $k=n-m+1$ then we have:
\[
    \Delta_{k}=(-1)^{k}b_{m}^{k}\sum_{j=0}^{n-m}t_{k-j}b_{(m-1)-j}
\]
Let $i=k-j,$ we see that
\[\begin{array}{lcl}
    \Delta_{k}&=&(-1)^{k}b_{m}^{k}\sum_{j=0}^{n-m}t_{k-j}b_{(m-1)-j}
        =(-1)^{k}b_{m}^{k}\sum_{i=k}^{1}t_{i}b_{(m-1)-(k-i)}\smallskip\\
        &=&(-1)^{k}b_{m}^{k}\sum_{i=1}^{k}t_{i}b_{m-(1+k)+i}
\end{array}
\]
as claim, and the theorem was proved.
\end{proof}

\begin{cor}\label{cor-3.4}
If $\Delta_{k}=\left|\begin{array}{lllll}
    b_{m-1}    &-b_{m}        &~~0          &\dots   &~~0\\
    b_{m-2}    &~~b_{m-1}     &-b_{m}       &\ddots  &~~\vdots\\
    \quad\vdots &\quad\vdots  &~~\ddots     &\ddots  &~~0\\
    b_{m-k+1}  &~~b_{m-k+2}   &             &\ddots  &-b_{m}\\
    b_{m-k}    &~~b_{m-k+1}   &~~b_{m-k+2}  &\dots   &~~b_{m-1}
    \end{array}\right|$ then
\begin{equation}\label{eq-Delta-k}
\Delta_{k}=b_{m}^{k}\displaystyle\sum_{i=1}^{k}t_{i}b_{m-(k+1)+i}
=b_{m}^{k}\displaystyle\sum_{j=1}^{k-1}t_{k-j}b_{(m-1)-j}
\end{equation}
where $\{t_{r}\}$ is the linear recurrent sequence,
$t_{1}=\frac{1}{b_{m}}$ and
$t_{r}=\frac{1}{b_{m}}\sum_{i=1}^{r-1}b_{m-i}t_{r-i}$\ for
$r=2,3,\ldots .$
\end{cor}

\begin{proof}
From Theorem \ref{thm-3.1} assert that
\begin{eqnarray}\label{eq56}
    \Delta_{k}&=&\left|\begin{array}{lllll}
    -b_{m-1}    &~~b_{m}    &~~0          &\dots   &~~0\\
    -b_{m-2}    &-b_{m-1}   &~~b_{m}      &\ddots  &~~\vdots\\
    \quad\vdots &\quad\vdots&~~\ddots     &\ddots  &~~0\\
    -b_{m-k+1}  &-b_{m-k+2} &             &\ddots  &~~b_{m}\\
    -b_{m-k}    &-b_{m-k+1} &-b_{m-k+2}   &\dots   &-b_{m-1}
    \end{array}\right|\nonumber\\
    &=&(-1)^{k}b_{m}^{k}\displaystyle\sum_{i=1}^{k}t_{i}b_{m-(k+1)+i}
    =(-1)^{k}b_{m}^{k}\displaystyle\sum_{j=1}^{k-1}t_{k-j}b_{(m-1)-j}
\end{eqnarray}
where $\{t_{r}\}$ is the linear recurrent sequence,
$t_{1}=\frac{1}{b_{m}}$ and
$t_{r}=\frac{1}{b_{m}}\sum_{i=1}^{r-1}b_{m-i}t_{r-i}$\ for
$r=2,3,\ldots.$ Multiplying (\ref{eq56}) both sides by
$(-1)^{k}$, we see that
\[
(-1)^{k}\left|\begin{array}{lllll}
    -b_{m-1}    &~~b_{m}    &~~0          &\dots   &~~0\\
    -b_{m-2}    &-b_{m-1}   &~~b_{m}      &\ddots  &~~\vdots\\
    \quad\vdots &\quad\vdots&~~\ddots     &\ddots  &~~0\\
    -b_{m-k+1}  &-b_{m-k+2} &             &\ddots  &~~b_{m}\\
    -b_{m-k}    &-b_{m-k+1} &-b_{m-k+2}   &\dots   &-b_{m-1}
    \end{array}\right|
\]
\[\begin{array}{lcl}
    \qquad &=&
\left|\begin{array}{lllll}
    b_{m-1}    &-b_{m}        &~~0          &\dots   &~~0\\
    b_{m-2}    &~~b_{m-1}     &-b_{m}       &\ddots  &~~\vdots\\
    \quad\vdots &\quad\vdots  &~~\ddots     &\ddots  &~~0\\
    b_{m-k+1}  &~~b_{m-k+2}   &             &\ddots  &-b_{m}\\
    b_{m-k}    &~~b_{m-k+1}   &~~b_{m-k+2}  &\dots   &~~b_{m-1}
    \end{array}\right|\smallskip\\
    \qquad&=&(-1)^{2k}b_{m}^{k}\displaystyle\sum_{i=1}^{k}t_{i}b_{m-(k+1)+i}
    =(-1)^{2k}b_{m}^{k}\displaystyle\sum_{j=1}^{k-1}t_{k-j}b_{(m-1)-j}\smallskip\\
    \qquad&=&b_{m}^{k}\displaystyle\sum_{i=1}^{k}t_{i}b_{m-(k+1)+i}
    =b_{m}^{k}\displaystyle\sum_{j=1}^{k-1}t_{k-j}b_{(m-1)-j}.
\end{array}
\]
The corollary was proved.
\end{proof}

\section*{Acknowledgments}

The author is very grateful to the anonymous referees for their comments and suggestions, which inspired
the improvement of the manuscript.

\bigskip

\end{document}